\documentclass[12pt]{article}
\usepackage{pgf,pgfarrows,pgfnodes,pgfautomata,pgfheaps,pgfshade}
\usepackage{amsmath,amssymb}
\usepackage{graphics}
\usepackage{amsmath}

\def\C{\mathbb{C}}
\newtheorem{theorem}{\hspace*{\parindent}Theorem}
\newtheorem{lemma}{\hspace*{\parindent}Lemma}
\newtheorem{corollary}{\hspace*{\parindent}Corollary}
\title{Hypergeometric functions as generalized Stieltjes transforms}
\author{D.\,Karp\footnote{Far Eastern Federal University, Vladivostok, Russia,
e-mail:\,\emph{dimkrp@gmail.com}}~~and~E.\,Prilepkina\footnote{Far
Eastern Federal University, Vladivostok, Russia,
e-mail:\,\emph{pril-elena@yandex.ru}}}
\date{}
\begin{document}
\maketitle

\begin{center}
\parbox{12cm}{
\small\textbf{Abstract.} In this paper we apply generalized
Stieltjes transform representation to study the generalized
hypergeometric function.  Among the results thus proved  are new
integral representations, inequalities, properties of the Pad\'{e}
table and the properties of the generalized hypergeometric
function as a conformal map. }
\end{center}

\bigskip


\bigskip


\bigskip

\paragraph{1. Introduction.}
Functions representable in one of the forms
\begin{equation}\label{eq:fmain}
f(z)=C_1+\int\limits_{[0,\infty)}\frac{\mu(du)}{(u+z)^{\alpha}}=\int\limits_{[0,\infty)}\frac{\rho(dt)}{(1+tz)^{\alpha}}+\frac{C_2}{z^{\alpha}},
\end{equation}
are known as generalized Stieltjes functions.  Here $\alpha>0$,
$\mu$ and $\rho$ are non-negative measures supported on
$[0,\infty)$, $C_1\geq{0}$, $C_2\geq{0}$ are constants and we
always choose the principal branch of the power function. The
measures $\mu$ and $\rho$ are assumed to produce convergent
integrals (\ref{eq:fmain}) for each
$z\in\C\!\setminus\!(-\infty,0]$ so that the function $f$ is
holomorphic in $\C\!\setminus\!(-\infty,0]$. Generalized Stieltjes
functions have been studied by a number of authors including
\cite{Sokal,Sumner}, \cite[Section~8]{WidderP},
\cite[Chapter~VIII]{Widder}. For more detailed overview of the
properties of generalized Stieltjes functions and related
bibliography see our recent paper \cite{KarpPril}.  In the same
paper we introduced the notion of the exact Stieltjes order as
follows.  If we define $S_{\alpha}$ to be the class of functions
representable by (\ref{eq:fmain}) then one can show that
$S_{\alpha}\subset{S_{\beta}}$ when $\alpha<\beta$.  We will say
that $f$ is of the exact Stieltjes order $\alpha^*$ if
$f\in\cup_{\alpha>0}S_{\alpha}$ and
\begin{equation}\label{eq:exactorder}
\alpha^*=\inf\{\alpha:~f\in{S_{\alpha}}\}.
\end{equation}
Using Sokal's characterization of $S_{\alpha}$ found in
\cite{Sokal} it is not difficult to see that $f\in{S_{\alpha^*}}$.
Moreover, in \cite{KarpPril} we gave a criterion of exactness
leading to some simple sufficient conditions.  In particular,  we
will need the following result contained in
\cite[Corollary~1]{KarpPril}.
\begin{theorem}\label{th:exactness}
Suppose $f\in{S_{\alpha}}$ and for sufficiently small
$\varepsilon>0$
\[
\lim\limits_{y\to+\infty}\frac{\Phi_{\varepsilon}(2y)}{\Phi_{\varepsilon}(y)}<1,
\]
where
\begin{equation}\label{eq:criterium}
\Phi_{\varepsilon}(y)=\int\limits_{(0,y)}\frac{\mu(du)}{(y-u)^{\varepsilon}}.
\end{equation}
Then $\alpha$ is the exact Stieltjes order of $f$.
\end{theorem}

In this paper we aim to apply the results of \cite{KarpPril} to
study the generalized hypergeometric function defined by the
series
\begin{equation}\label{eq:pFqdefined}
{_{q+1}F_q}\left(\left.\!\!\begin{array}{c}\sigma, A\\
B\end{array}\right|z\!\right)={_{q+1}F_q}\left(\sigma,A;B;z\right):=\sum\limits_{n=0}^{\infty}\frac{(\sigma)_n(a_1)_n(a_2)_n\cdots(a_{q})_n}{(b_1)_n(b_2)_n\cdots(b_q)_nn!}z^n,
\end{equation}
where we write $A=(a_1,a_2,\ldots,a_q)$, $B=(b_1,b_2,\ldots,b_q)$
for brevity and $(a)_0=1$, $(a)_n=a(a+1)\cdots(a+n-1)$,
$n\geq{1}$,  denotes the rising factorial. The series
(\ref{eq:pFqdefined}) converges in the unit disk and its sum can
be extended analytically  to the whole complex plane cut along the
ray $[1,\infty)$. See details in \cite{AAR,Kiryakova,PBM3}.

Euler's integral representation \cite[Theorem~2.2.1]{AAR}
$$
{_2F_1}(\sigma,a;b;-z)\!=\!\frac{\Gamma(b)}{\Gamma(a)\Gamma(b-a)}\!\int\limits_{0}^{1}\frac{t^{a-1}(1-t)^{b-a-1}dt}{(1+zt)^{\sigma}}
\!=\!\frac{\Gamma(b)}{\Gamma(a)\Gamma(b-a)}\!\int\limits_{1}^{\infty}\frac{u^{\sigma-b}(u-1)^{b-a-1}du}{(u+z)^{\sigma}}
$$
for the Gauss hypergeometric function $_2F_1$ shows that it is a
generalized Stieltjes function at least when $b>a>0$ and
$\sigma>0$. In her book \cite{Kiryakova} Virginia Kiryakova gave
the representation
$$
{_{q+1}F_q}\left(\left.\!\!\begin{array}{c}\sigma, A\\
B\end{array}\right|-z\!\right)
=\!\!\int\limits_{0}^{1}\frac{\rho(s)ds}{(1+sz)^{\sigma}}
$$
under the constraints $b_k>a_k>0$, $k=1,2,\ldots,q$, and with
$\rho$ expressed in terms of Meijer's $G$-function (see
(\ref{eq:rho-G}) below).  In \cite{KarpSitnik} Karp and Sitnik
established  the same formula but with $\rho$ expressed by a
multidimensional integral which is manifestly positive under the
same constraints. In this work we generalize both these results by
stating necessary and sufficient conditions for the above
representation to hold and sufficient conditions for the weight
$\rho$ to be non-negative (the latter conditions are also believed
to be necessary but we have no proof of this claim). We find the
exact Stieltjes order of ${_{q+1}F_q}$ and give a number of
consequences, including new integral representations,
inequalities, properties of the Pad\'{e} table and properties of
${_{q+1}F_q}$ as a conformal map.

\paragraph{2. The exact Stieltjes order of $_{q+1}F_q$.}

We will  need a particular case of the Meijer's $G$-function
defined by (see \cite{Kiryakova,PBM3})
\begin{equation}\label{eq:G-defined}
G^{q,0}_{p,q}\!\left(\!z~\vline\begin{array}{l}a_1,\ldots,a_p
\\b_1,\ldots, b_q\end{array}\!\!\right):=
\frac{1}{2\pi{i}}
\int\limits_{c-i\infty}^{c+i\infty}\!\!\frac{\Gamma(b_1\!+\!s)\dots\Gamma(b_q\!+\!s)}
{\Gamma(a_{1}\!+\!s)\dots\Gamma(a_p\!+\!s)}z^{-s}ds,
\end{equation}
where $c>-\min(\Re{b_1},\Re{b_2},\ldots,\Re{b_q})$.  Since the
gamma function is real symmetric,
$\Gamma(\overline{z})=\overline{\Gamma(z)}$, the function
$G^{q,0}_{p,q}$ is  real if all parameters $a_i$, $b_i$ are real.
Define
\begin{equation}\label{eq:psi-defined}
\psi:=\sum\limits_{k=1}^{q}(b_k-a_k).
\end{equation}

\begin{lemma}\label{lm:G-zero}
Set $A=(a_1,\ldots,a_q)$, $B=(b_1,\ldots,b_q)$. If
\begin{equation}\label{eq:Reba-cond}
\Re(\psi)>0,
\end{equation}
then
\begin{equation}\label{eq:G-zero}
G^{q,0}_{q,q}\left(x\left|\begin{array}{l}\!\!B\!\!\\\!\!A\!\!\end{array}\right.\right)=0
~~\text{for}~~~x>1.
\end{equation}
\end{lemma}
\textbf{Proof.} From (\ref{eq:G-defined}) we have
$$
G^{q,0}_{q,q}\left(x\left|\begin{array}{l}\!\!B\!\!\\\!\!A\!\!\end{array}\right.\right)
=\frac{1}{2\pi{i}}\lim\limits_{R\to\infty}
\int\limits_{c-iR}^{c+iR}
\frac{\Gamma(a_1+s)\cdots\Gamma(a_q+s)}{\Gamma(b_1+s)\cdots\Gamma(b_q+s)}
e^{-s\ln{x}}ds.
$$
Expression under the integral sign has no poles inside the closed
contour starting at the point $c-iR$, tracing the semicircle
$c+Re^{i\varphi}$, $-\pi/2\leq\varphi\leq\pi/2$, upto the point
$c+iR$ and then back to $c-iR$ along the line segment $c+it$,
$-R\leq{t}\leq{R}$. Hence, we have by the Cauchy theorem:
\[
I(R):=\frac{1}{2\pi}\int\limits_{-R}^{R}
\frac{\Gamma(a_1+c+it)\cdots\Gamma(a_q+c+it)}{\Gamma(b_1+c+it)\cdots\Gamma(b_q+c+it)}
e^{-(c+it)\ln{x}}dt=
\]
\[
=-\frac{R}{2\pi}e^{-c\ln{x}}\int\limits_{-\pi/2}^{\pi/2}
\frac{\Gamma(a_1+c+Re^{i\varphi})\cdots\Gamma(a_q+c+Re^{i\varphi})}
{\Gamma(b_1+c+Re^{i\varphi})\cdots\Gamma(b_q+c+Re^{i\varphi})}e^{i(R\ln{x}\sin{\varphi}+\varphi)}
e^{-R\ln{x}\cos{\varphi}}d\varphi.
\]
Set $z=Re^{i\varphi}$. Using Stirling's asymptotic formula (see,
for instance, \cite[Theorem~1.4.2]{AAR}) we get the relation
\[
\log\left\{\frac{\Gamma(a_1+c+z)\cdots\Gamma(a_q+c+z)}{\Gamma(b_1+c+z)\cdots\Gamma(b_q+c+z)}\right\}
=-\psi\log(z)+O(1/z)~\text{as}~|z|\to\infty,
\]
which holds uniformly in the sector $|\arg{z}|\leq\pi-\delta$, for
each $\delta\in(0,\pi)$. Hence,
\[
\left|\frac{\Gamma(a_1+c+z)\cdots\Gamma(a_q+c+z)}{\Gamma(b_1+c+z)\cdots\Gamma(b_q+c+z)}\right|
=R^{-\Re(\psi)}(1+O(1/R)),~~~R\to\infty.
\]
Consequently,
$$
|I(R)|=O\left(R^{-\Re(\psi)+1}\right)\int\limits_{-\pi/2}^{\pi/2}
e^{-R\ln{x}\cos{\varphi}}d\varphi ~~~\text{as}~~~R\to\infty.
$$
Applying the inequality $\cos\varphi\geq 1-\frac{2}{\pi}\varphi$,
$0\leq\varphi\leq\pi/2$, we obtain (recall that $x>1$)
\[
\int\limits_{0}^{\pi/2}e^{-R\ln{x}\cos{\varphi}}d\varphi\leq
\int\limits_{0}^{\pi/2}e^{-R\left(1-\frac{2}{\pi}\varphi\right)\ln{x}}
d\varphi=e^{-R\ln{x}}\int\limits_{0}^{\pi/2}e^{\frac{2}{\pi}R\varphi\ln{x}}
d\varphi=\frac{\pi}{2R\ln{x}}\left(1-e^{-R\ln{x}}\right).
\]
Combining this estimate with the previous relation  we see that
\[
\lim\limits_{R\to\infty}I(R)=0~~~\text{for each}~x>1.~~~\square
\]

\textbf{Remark.} Formula (\ref{eq:G-zero}) is given in
\cite[formula (8.2.2.2)]{PBM3} under more restrictive conditions
then (\ref{eq:Reba-cond}).  For this reason we decided to include
a direct proof.

\begin{theorem}\label{th:FGrep}
Suppose $|\arg(1+z)|<\pi$ and $\sigma$ is  an arbitrary complex
number. Representation
\begin{equation}\label{eq:Frepr}
{_{q+1}F_q}\left(\left.\!\!\begin{array}{c}\sigma,A\\
B\end{array}\right|-z\!\right)
=\!\!\int\limits_{0}^{1}\frac{\rho(s)ds}{(1+sz)^{\sigma}}
\end{equation}
with a summable on $[0,1]$ function $\rho$ holds true if and only
if $\Re{a_i}>0$ for $i=1,\ldots,q$ and $\Re{\psi}>0$, where $\psi$
is defined in \emph{(\ref{eq:psi-defined})}. Under these
conditions
\begin{equation}\label{eq:rho-G}
\rho(s)=\left(\prod\limits_{i=1}^{q}\frac{\Gamma(b_i)}{\Gamma(a_i)}\right)\frac{1}{s}G^{q,0}_{q,q}\left(s\left|\begin{array}{l}\!\!B\!\!\\\!\!A\!\!\end{array}\right.\right).
\end{equation}
\end{theorem}
\textbf{Remark.} Representation (\ref{eq:Frepr}) after change of
variable $t=1/s$ can also be written as
\begin{equation}\label{eq:Frepr1}
{_{q+1}F_q}\left(\left.\!\!\begin{array}{c}\sigma,A\\
B\end{array}\right|-z\!\right)
=\!\!\int\limits_{1}^{\infty}\frac{\mu(t)dt}{(t+z)^{\sigma}},
\end{equation}
\begin{equation}\label{eq:mu-G}
\mu(t)=\left(\prod\limits_{i=1}^{q}\frac{\Gamma(b_i)}{\Gamma(a_i)}\right)t^{\sigma-1}
G^{q,0}_{q,q}\left(1/t\left|\begin{array}{l}\!\!B\!\!\\\!\!A\!\!\end{array}\right.\right)
\end{equation}
- a form which we will also use.

\noindent\textbf{Proof.} Suppose first that $\Re{a_i}>0$ for
$i=1,\ldots,q$ and $\Re{\psi}>0$. Consider the right-hand side of
(\ref{eq:Frepr}) with $\rho$ given by (\ref{eq:rho-G}). Applying
the binomial expansion to $(1+sz)^{-\sigma}$  and integrating term
by term we immediately obtain the left-hand side of
(\ref{eq:Frepr}) since
\[
\int\limits_{0}^{1}s^{k}\rho(s)ds=\int\limits_{0}^{\infty}s^{k}\rho(s)ds=\frac{(a_1)_k\cdots(a_{q})_k}{(b_1)_k\cdots(b_q)_k}.
\]
The first equality here is due to Lemma~\ref{lm:G-zero}.  The
second equality expresses the basic property of the Meijer's
$G$-function: its Mellin transform is equal to the ratio of the
appropriate gamma functions (see, for instance, \cite[formula
2.24.2.1]{PBM3} or \cite[formula (A.25), p.319]{Kiryakova}). The
integral converges uniformly in $k$ in the neighbourhood of $s=0$
since
\begin{equation}\label{G-asymp-zero}
G^{q,0}_{q,q}\!\left(s~\vline\begin{array}{l}\!B\!\\\!A\!\end{array}\right)=
O\left(s^{a}\ln^{m-1}(1/s)\right),~~s\to{0},
\end{equation}
where $a=\min(\Re(a_1),\ldots,\Re(a_q))>0$ by assumption and the
minimum is taken over those $a_i$ for which there is no
$b_j=a_i-l$ for some $l\in\mathbb{N}_0$. The minimum can be
attained for several different numbers $a_i$ and then $m$ is the
maximal multiplicity among these numbers. This formula follows
from \cite[Corollary~1.12.1]{KilSaig} or \cite[formula
(11)]{Karp}. The integral converges uniformly in $k$ in the
neighbourhood of $s=1$ because, the function $G^{q,0}_{q,q}$ has a
singularity of the magnitude $(1-s)^{\Re(\psi)-1}$ possibly
multiplied by logarithmic terms if $\Re(\psi)\leq{1}$ and is
bounded if $\Re(\psi)>1$ (see \cite[8.2.59]{PBM3}). Hence,
condition (\ref{eq:Reba-cond}) guarantees   uniform integrability
of $\rho$ in the neighbourhood of $s=1$.  Uniform integrability
justifies the interchange of summation and integration.

To prove necessity suppose that (\ref{eq:Frepr}) holds with a
summable function $\rho$. Then
\begin{equation}\label{eq:rho-moments}
\int\limits_{0}^{1}s^{k}\rho(s)ds=\frac{(a_1)_k\cdots(a_{q})_k}{(b_1)_k\cdots(b_q)_k}
\end{equation}
by termwise integration and comparing with (\ref{eq:pFqdefined}).
We aim to show that $\Re{a_i}>0$ for $i=1,\ldots,q$ and
$\Re{\psi}>0$. Assume first that $\Re{a_i}\leq{0}$ for some $i$
while $\Re{\psi}>0$. The asymptotic formula (\ref{G-asymp-zero})
combined with Lemma~\ref{lm:G-zero} shows that
$$
\left(\prod\limits_{i=1}^{q}\frac{\Gamma(b_i)}{\Gamma(a_i)}\right)\int\limits_{0}^{1}s^{k-1}G^{q,0}_{q,q}\!\left(s~\vline\begin{array}{l}\!B\!\\\!A\!\end{array}\right)ds
=\frac{(a_1)_k\cdots(a_{q})_k}{(b_1)_k\cdots(b_q)_k}
$$
for $k>-a$, where as before $a=\min(\Re(a_1),\ldots,\Re(a_q))$.
Hence all moments of the functions $s^{[-a]+1}\rho$ and
$s^{[-a]}G^{q,0}_{q,q}$ coincide. This implies that $\rho$ must be
given by (\ref{eq:rho-G}) by the determinacy of the moment problem
on a finite interval.  But then the integral in (\ref{eq:Frepr})
must diverge by (\ref{G-asymp-zero}).  A contradiction.

If $\Re{\psi}<0$ the sequence
$$
\frac{(a_1)_k\cdots(a_{q})_k}{(b_1)_k\cdots(b_q)_k}
$$
is unbounded and cannot serve as a moment sequence of a signed
measure on $[0,1]$, so that (\ref{eq:rho-moments}) is impossible
and hence so is (\ref{eq:Frepr}). Finally, if $\Re{\psi}=0$ a
careful application of Stirling's formula shows that this sequence
tends to a non-zero constant as $k\to\infty$ (see \cite[formula
(1.2.5)]{KilSaig}) while the left-hand side of
(\ref{eq:rho-moments}) must tend to zero for any summable function
$\rho$, so again a contradiction.~~$\square$

\textbf{Remark.} Formula (\ref{eq:Frepr}) has been discovered by
Kiryakova in \cite{Kiryakova} by iterative fractional integrations
under additional assumption that all parameters are real and
$b_k>a_k>0$, $k=1,2,\ldots,q$. The elementary proof included here
is not contained in this reference.

\medskip

In the sequel we will need the notion of majorization
\cite[Definition A.2, formula (12)]{MOA}. It is said that
$B=(b_1,\ldots,b_q)$ is weakly supermajorized by
$A=(a_1,\ldots,a_q)$ (symbolized by $B\prec^W{A}$) if
\begin{equation}\label{eq:amajorb}
\begin{split}
& 0<a_1\leq{a_2}\leq\cdots\leq{a_q},~~
0<b_1\leq{b_2}\leq\cdots\leq{b_q},
\\
&\sum\limits_{i=1}^{k}a_i\leq\sum\limits_{i=1}^{k}b_i~~\text{for}~~k=1,2\ldots,q.
\end{split}
\end{equation}
If in addition $\psi(=\sum_{i=1}^{q}(b_i-a_i))=0$ then $B$ is said
to be majorized by $A$, $B\prec{A}$.

\begin{lemma}\label{lm:G-positive}
Suppose that $B\prec^W{A}$ but not $B\prec{A}$ \emph{(}that is
$\psi>0$\emph{)}. Then for all $0<s<1$
\begin{equation}\label{eq:G-positive}
G^{q,0}_{q,q}\left(s~\vline\begin{array}{l}\!B\\\!A\end{array}\!\right)\geq{0}.
\end{equation}
\end{lemma}
\textbf{Proof.}  Alzer showed in \cite[Theorem~10]{Alzer} that the
function
\[
x\to\prod\limits_{i=1}^{q}\frac{\Gamma(x+a_i)}{\Gamma(x+b_i)}
\]
is completely monotonic on $(0,\infty)$ if $B\prec^W{A}$. This
implies that the sequence
\[
\left\{\prod\limits_{i=1}^{q}\frac{\Gamma(n+a_i)}{\Gamma(n+b_i)}\right\},~~~n=0,1,2,\ldots,
\]
is a completely monotonic sequence.  Hence by the Hausdorff
theorem there exists a unique non-negative measure $d\nu$
supported on $[0,1]$ such that
\[
\int\limits_{[0,1]}s^nd\nu(s)=\prod\limits_{i=1}^{q}\frac{\Gamma(n+a_i)}{\Gamma(n+b_i)}.
\]
On the other hand if $\psi>0$
\[
\int\limits_{0}^{1}s^{n-1}G^{q,0}_{q,q}\left(s\,\,\vline\begin{array}{c}\!B\\
\!A\end{array}\!\!\right)ds=\prod\limits_{i=1}^{q}\frac{\Gamma(n+a_i)}{\Gamma(n+b_i)},
\]
so that by determinacy of the Haudorff moment problem
\[
d\nu(s)=\frac{1}{s}G^{q,0}_{q,q}\left(s\,\,\vline\begin{array}{c}\!B\\
\!A\end{array}\!\!\right)ds.
\]
Non-negativity of the measure completes the proof.~$\square$

\textbf{Remark.}  According to Bernstein's theorem every
completely monotonic function on $(0,\infty)$ is the Laplace
transform of a non-negative measure.  The proof of
Lemma~\ref{lm:G-positive} shows that the representing measure in
Alzer's theorem~10 from \cite{Alzer} is given by
$$
\prod\limits_{i=1}^{q}\frac{\Gamma(x+a_i)}{\Gamma(x+b_i)}=\int\limits_{0}^{\infty}e^{-tx}
G^{q,0}_{q,q}\left(e^{-t}\,\,\vline\begin{array}{c}\!B\\
\!A\end{array}\!\!\right)\!dt.
$$

\textbf{Remark.} By taking the Mellin transform on both sides and
changing variables one can show that for $x>0$
\begin{multline}\label{eq:G-multiint}
G^{q,0}_{q,q}\left(x\,\,\vline\begin{array}{c}\!B\\
\!A\end{array}\!\!\right)
=\frac{x^{a_1}}{\prod_{i=1}^{q}\Gamma(b_i-a_i)}\\
\times\int\limits_{\Lambda_q(x)}[1-x/(t_2\cdots{t_q})]^{b_1-a_1-1}\prod_{k=2}^{q}t_k^{a_k-a_1-1}(1-t_k)^{b_k-a_k-1}\,dt_2\cdots{dt_q},
\end{multline}
if $\Re(b_k-a_k)>0$, $k=1,2,\ldots,q$, $q\geq{2}$. Here the domain
of integration is given by
\begin{equation}\label{eq:Lambda}
\Lambda_q(x)=[0,1]^{q-1}\cap\{t_2,\ldots,t_q:~t_2\cdots{t_q}>x\}.
\end{equation}
This formula shows the positivity of $G^{q,0}_{q,q}$ under the
conditions $b_k>a_k>0$, $k=1,2,\ldots,q$, which are manifestly
more restrictive then $B\prec^W{A}$ and $\psi>0$. Formula
(\ref{eq:G-multiint}) is implicit in \cite{KarpSitnik}.

\begin{theorem}\label{th:q1Fq-order}
Suppose $0<\sigma\leq\min(a_1,\ldots,a_q)$ and $B\prec^W{A}$. Then
$f:={_{q+1}F_q}(\,\sigma,A;B;-z)$ is a generalized Stieltjes
function of the exact order $\sigma$.  In particular, $f$ is
completely monotonic.
\end{theorem}
\textbf{Proof.} Assume first that
$\psi(=\sum_{i=1}^{q}(b_i-a_i))>0$. Then by Theorem~\ref{th:FGrep}
$f$ is represented by (\ref{eq:Frepr1}) with the measure $\mu$
non-negative by Lemma~\ref{lm:G-positive}. Hence,
$f\in{S_{\sigma}}$.  To show that $\sigma$ is exact we will apply
Theorem~\ref{th:exactness}. Fixing $\varepsilon>0$ compute
\[
\Phi_{\epsilon}(y):=\int\limits_{1}^{y}\frac{
\mu(u)du}{(y-u)^{\varepsilon}},
\]
where $\mu(u)$ is given by (\ref{eq:mu-G}).  Changing variable
$\tau=1/u$ and manipulating a little we obtain
\[
\Phi_{\epsilon}(y)=\frac{1}{y^\varepsilon}\int\limits_{1/y}^{1}
\frac{\tau^{\varepsilon-1-\sigma}}{(\tau-1/y)^\varepsilon}G^{q,0}_{q,q}\left(\tau\,\,\vline\begin{array}{c}\!B\\
\!A\end{array}\!\!\right)\!d\tau.
\]
According to \cite[formula (2.24.3)]{PBM3}  combined with
(\ref{eq:G-zero}) we get
\begin{equation}\label{eq:FyviaG}
\Phi_{\epsilon}(y)=\frac{\Gamma(1-\varepsilon)}{y^{\varepsilon}}
y^{\sigma}G^{q+1,0}_{q+1,q+1}\left(\frac{1}{y}\,\,\vline\begin{array}{c}1-\varepsilon+\sigma, B\\
\sigma, A\end{array}\!\!\right).
\end{equation}
Using (\ref{G-asymp-zero}) for the main asymptotic term of
$G^{q+1,0}_{q+1,q+1}$ we immediately arrive at
\begin{equation}\label{eq:limFy}
\lim\limits_{y\to+\infty}\frac{\Phi_{\varepsilon}(2y)}{\Phi_{\varepsilon}(y)}=2^{-\varepsilon}<1.
\end{equation}
Hence, by Theorem~\ref{th:exactness} the order $\sigma$ is exact.

Next, suppose that $B\prec{A}$, i.e. (\ref{eq:amajorb}) holds with
$\psi=0$. By Alzer's theorem the sequence on the right of
(\ref{eq:rho-moments}) is still a moment sequence of a
non-negative measure (see proof of Lemma~\ref{lm:G-positive})
which shows that $f\in{S_{\sigma}}$. We will, however, give
another proof of this fact which will extend to a proof of the
exactness of $\sigma$. Consider the sequence
\[
f_m(z)={_{q+1}F_q}(\,\sigma,A;B',b_q+1/m;-z),~~~B'=(b_1,\ldots,b_{q-1}).
\]
According to what we have just proved each $f_m\in{S_{\sigma}}$
and the order is exact. We aim to apply
\cite[Theorem~10]{KarpPril} to show that $f\in{S_{\sigma}}$. To
this end we need to demonstrate that $f_m(x)\to{f(x)}$ for all
$x>0$. If $|z|<1$ then
\[
|f(z)-f_m(z)|\leq\sum\limits_{k=0}^{\infty}\frac{(\sigma)_k(a_1)_k\cdots(a_q)_k|z|^k}{(b_1)_k\cdots(b_{q-1})_kk!}\left[\frac{1}{(b_q)_k}-\frac{1}{(b_q+1/m)_k}\right]
\underset{m\to\infty}{\to}0
\]
due to uniform  in $m$ convergence of the series.  The convergence
can be extended to all $z\in\C\!\setminus(-\infty,-1]$ using
Vitali-Porter (or Stieltjes-Vitali) theorem on induced convergence
\cite[Corollary~7.5]{Burckel}.  This theorem requires the set
$\{f_m\}$ to be locally uniformly bounded in
$\C\!\setminus(-\infty,-1]$. This boundedness can be seen from the
easily verifiable contiguous relation
\begin{multline}
{_{q+1}F_q}(\,\sigma,A;B',b_q+1/m;-z)={_{q+1}F_q}(\,\sigma,A;B',b_q+1+1/m;-z)
\\
-\frac{z\sigma\prod_{i=1}^{q}a_i}{(b_q+1/m)(b_q+1+1/m)\prod_{i=1}^{q-1}b_i}{_{q+1}F_q}(\,\sigma,A+1;B'+1,b_q+2+1/m;-z),
\end{multline}
where both functions on the right are bounded uniformly in $m$ due
to representation (\ref{eq:Frepr}).  This proves that
$f\in{S_{\sigma}}$. Finally, we need to demonstrate that the order
$\sigma$ is exact for $f$.  The distribution function of the
representing measure of $f_m$ is given by
\begin{multline*}
F_m(y)=\prod\limits_{i=1}^{q}\frac{\Gamma(b_i)}{\Gamma(a_i)}\int\limits_{[1,y)}t^{\sigma-1}
\mathrm{G}^{q,0}_{q,q}\left(1/t\,\,\vline\begin{array}{c}B',
b_q+1/m\\A\end{array}\!\!\right)dt
\\
=\prod\limits_{i=1}^{q}\frac{\Gamma(b_i)}{\Gamma(a_i)}y^{\sigma}\mathrm{G}^{q+1,0}_{q+1,q+1}\left(1/y\,\,\vline\begin{array}{c}1+\sigma,B',
b_q+1/m\\\sigma, A\end{array}\!\!\right),
\end{multline*}
where we again used \cite[formula (2.24.3)]{PBM3} combined with
(\ref{eq:G-zero}).  Taking limit as $m\to\infty$ we obtain the
distribution function of the measure representing $f$ in the form
$$
\prod\limits_{i=1}^{q}\frac{\Gamma(b_i)}{\Gamma(a_i)}y^{\sigma}\mathrm{G}^{q+1,0}_{q+1,q+1}\left(1/y\,\,\vline\begin{array}{c}1+\sigma,B\\\sigma,
A\end{array}\!\!\right).
$$
Comparing this formula with (\ref{eq:FyviaG}) for $\varepsilon=0$
we see that the distribution function does not change its form
whether $\psi>0$ or $\psi=0$.  This implies that the function
$\Phi_{\varepsilon}(y)$ is again expressed by (\ref{eq:FyviaG})
when $\psi=0$ (since $\Phi_{\varepsilon}$ is proportional to the
fractional derivative of order $\varepsilon$ of the distribution
function).  Hence, the limit in (\ref{eq:limFy}) is again less
than 1 which according to Theorem~\ref{th:exactness} proves the
exactness of the order $\sigma$.~~ $\square$

\textbf{Remark.} If $\psi>0$ then the representing measure in the
above theorem is given in (\ref{eq:Frepr}) or (\ref{eq:Frepr1}).
However, if $B\prec{A}$ (i.e. $\psi=0$) then
Theorem~\ref{th:q1Fq-order} leaves the question of finding the
representing measure open.  For $q=1$ the answer is obvious:
$$
{_{2}F_1}(\,\sigma,a;a;-z)=\frac{1}{(1+z)^{\sigma}}
$$
by the  binomial theorem, so that the representing measure is
$\delta_1$ (the Dirac measure concentrated at $1$).
 For $q=2$ representation (\ref{eq:Frepr}) reduces to (see \cite[Lemma~2]{KarpSitnik})
 \begin{multline}\label{eq:3F2-int2F1}
{_3F_2}(\sigma,a_1,a_2;b_1,b_2;-z)=\frac{\Gamma(b_1)\Gamma(b_2)}{\Gamma(a_1)\Gamma(a_2)\Gamma(b_1+b_2-a_1-a_2)}
\\
\times\int\limits_{0}^{1}\frac{t^{a_2-1}(1-t)^{b_1+b_2-a_1-a_2-1}}{(1+zt)^{\sigma}}{_2F_1}(b_1-a_1,b_2-a_1;b_1+b_2-a_1-a_2;1-t)dt
\end{multline}
valid if $a_1,a_2>0$, $b_1+b_2>a_1+a_2$. To compute the limiting
measure when $b_1+b_2=a_1+a_2$ we put $\epsilon=b_1+b_2-a_1-a_2$,
$\varphi(t)=(1+zt)^{\sigma}$ and let $\epsilon\to{0}$ in
\begin{multline*}
\frac{1}{\Gamma(\epsilon)}
\int\limits_{0}^{1}t^{a_2-1}(1-t)^{\epsilon-1}{_2F_1}(b_1-a_1,b_2-a_1;\epsilon;1-t)\varphi(t)dt
\\
=\frac{1}{\Gamma(\epsilon)}
\int\limits_{0}^{1}(1-u)^{a_2-1}u^{\epsilon-1}{_2F_1}(b_1-a_1,b_2-a_1;\epsilon;u)\psi(u)du,
\end{multline*}
where $t=1-u$ and $\psi(u):=\varphi(1-u)$.  We have
$$
{_2F_1}(b_2-a_1,b_1-a_1;\epsilon;u)
=1+\frac{\Gamma(\epsilon)}{\Gamma(b_2-a_1)\Gamma(b_1-a_1)}
\sum\limits_{k=1}^{\infty}\frac{\Gamma(b_2-a_1+k)\Gamma(b_1-a_1+k)}{\Gamma(\epsilon+k)k!}u^k.
$$
Hence,
\begin{multline*}
\lim\limits_{\epsilon\to{0}}\frac{1}{\Gamma(\epsilon)}
\int\limits_{0}^{1}(1-u)^{a_2-1}u^{\epsilon-1}{_2F_1}(b_1-a_1,b_2-a_1;\epsilon;u)\psi(u)du
\\
=\lim\limits_{\epsilon\to{0}}\!\frac{1}{\Gamma(\epsilon)}
\!\int\limits_{0}^{1}\!(1-u)^{a_2-1}u^{\epsilon-1}\!\left[1\!+\!\frac{\Gamma(\epsilon)}{\Gamma(b_2-a_1)\Gamma(b_1-a_1)}
\!\sum\limits_{k=1}^{\infty}\frac{\Gamma(b_2-a_1+k)\Gamma(b_1-a_1+k)}{\Gamma(\epsilon+k)k!}u^k\right]\!\psi(u)du
\\
=\lim\limits_{\epsilon\to{0}}\frac{1}{\Gamma(\epsilon)}
\int\limits_{0}^{1}(1-u)^{a_2-1}u^{\epsilon-1}(\psi(0)+u\psi'(0)+O(u^2))du
\\
+\frac{1}{\Gamma(b_2-a_1)\Gamma(b_1-a_1)}\lim\limits_{\epsilon\to{0}}
\int\limits_{0}^{1}(1-u)^{a_2-1}u^{\epsilon}
\left[\sum\limits_{k=1}^{\infty}\frac{\Gamma(b_2-a_1+k)\Gamma(b_1-a_1+k)}{\Gamma(\epsilon+k)k!}u^{k-1}\right]\psi(u)du
\\
=\psi(0)\lim\limits_{\epsilon\to{0}}\frac{\Gamma(a_2)\Gamma(\epsilon)}{\Gamma(a_2+\epsilon)\Gamma(\epsilon)}+
\psi'(0)\lim\limits_{\epsilon\to{0}}\frac{\Gamma(a_2)\Gamma(\epsilon+1)}{\Gamma(a_2+\epsilon+1)\Gamma(\epsilon)}+\cdots
\\
+\frac{1}{\Gamma(b_2-a_1)\Gamma(b_1-a_1)}
\int\limits_{0}^{1}(1-u)^{a_2-1}
\left[\sum\limits_{k=1}^{\infty}\frac{\Gamma(b_2-a_1+k)\Gamma(b_1-a_1+k)}{(k-1)!k!}u^{k-1}\right]\psi(u)du
\\
=\psi(0)+ (b_2-a_2)(b_1-a_1)\int\limits_{0}^{1}(1-u)^{a_2-1}
\left[\sum\limits_{k=0}^{\infty}\frac{(b_2-a_1+1)_k(b_1-a_1+1)_k}{(2)_kk!}u^{k}\right]\psi(u)du.
\end{multline*}
Summing the series we get
$$
\psi(0)+ (b_1-a_1)(b_2-a_2)\int\limits_{0}^{1}(1-u)^{a_2-1}
{_2F_1}(b_1-a_1+1,b_2-a_1+1;2;u)\psi(u)du.
$$
So we have the following result: if $b_1+b_2=a_1+a_2$ then the
representing measure has an atom at $t=1$ ($\psi(0)=\varphi(1)$)
and a continuous part given above, so that
\begin{multline*}
{_{3}F_2}(\sigma,a_1,a_2;b_1,b_2,-z)=\frac{\Gamma(b_1)\Gamma(b_2)}{\Gamma(a_1)\Gamma(a_2)}\biggl\{\frac{1}{(1+z)^{\sigma}}
\\
+\int\limits_{0}^{1}\frac{(b_2-a_1)(b_1-a_1)t^{a_2-1}}{(1+zt)^{\sigma}}{_{2}F_1}(b_1-a_1+1,b_2-a_1+1;2;1-t)dt\biggl\}.
\end{multline*}
This formula can also be proved by comparing power series
coefficients on both sides and using the Gauss summation theorem.
Finding the representing measure for general $q$ remains an
interesting open problem we plan to deal with in a separate
publication.

\begin{corollary}\label{cr:q1Fqnewrepr}
Suppose $B\prec^W{A}$ with $\psi>0$, $\sigma\geq{2}$ and
$|\arg(z)|<\pi/\sigma$. Then
\begin{equation}\label{eq:q1Fqnewrepr}
\begin{split}
&{_{q+1}F_q}(\sigma,A,B;-z)=\int\limits_{0}^{\infty}\frac{\varphi(y)dy}{y^{\sigma}+z^{\sigma}},~\text{where}
\\
&\varphi(y)=\frac{\sigma{y^{\sigma-1}}}{\pi}\left(\prod\limits_{i=1}^{q}\frac{\Gamma(b_i)}{\Gamma(a_i)}\right)\int\limits_{0}^{1}
\frac{\sin\left\{\sigma\arctan\left(\frac{ty\sin(\pi/\sigma)}{1+ty\cos(\pi/\sigma)}\right)\right\}}
{t\left(1+2ty\cos(\pi/\sigma)+t^2y^2\right)^{\sigma/2}}G^{q,0}_{q,q}\left(t\,\,\vline\begin{array}{c}\!B\\
\!A\end{array}\!\!\right)dt.
\end{split}
\end{equation}
\end{corollary}
\textbf{Proof.}  According to \cite[Theorem~13]{KarpPril} combined
with Theorem~\ref{th:q1Fq-order} above the function\\
${_{q+1}F_q}(\sigma,A;B;-z^{1/\sigma})$ belongs to $S_1$ for
$\sigma>1$ under the assumptions of the corollary.  According to
the Stieltjes inversion formula \cite[Chapter VIII, Theorem
7b]{Widder}  the density of the representing measure for
$f\in{S_1}$ is found from ($x>0$):
\[
\frac{1}{2{\pi}i}\lim\limits_{\varepsilon\to{0}}[f(-x-i\varepsilon)-f(-x+i\varepsilon)].
\]
Substituting the the first formula (\ref{eq:Frepr}) for $f$ and
computing the limit we arrive at
(\ref{eq:q1Fqnewrepr}).~~$\square$

\textbf{Remark.}  For $1<\sigma<2$ a similar formula can be
obtained. However, since it's more cumbersome than
(\ref{eq:q1Fqnewrepr}) we decided to omit it.

\textbf{Remark.} Using the identity $\sin(2\arctan(s))=2s/(1+s^2)$
formula (\ref{eq:q1Fqnewrepr}) for $\sigma=2$  simplifies to
($|\arg(z)|<\pi/2$)
\begin{equation}
\begin{split}
&{_{q+1}F_q}(2,A,B;-z)=\int\limits_{0}^{\infty}\frac{\varphi(y)dy}{y^{2}+z^{2}},~~\text{where}
\\
&\varphi(y)=\frac{4}{\pi}\left(\prod\limits_{i=1}^{q}\frac{\Gamma(b_i)}{\Gamma(a_i)}\right)\int\limits_{0}^{1}
\frac{y^2}
{\left(1+t^2y^2\right)^{2}}G^{q,0}_{q,q}\left(t\,\,\vline\begin{array}{c}\!B\\
\!A\end{array}\!\!\right)\!dt.
\end{split}
\end{equation}

\textbf{Remark.} For the Gauss hypergeometric function formula
(\ref{eq:q1Fqnewrepr}) reduces to
\[
{_{2}F_1}(a,b;c;-z)=\int\limits_{0}^{\infty}\frac{\varphi(y)dy}{y^a+z^a},~~~c>b>0,~~a\geq{2}.
\]\[
\varphi(y)=\frac{a\Gamma(c)y^{a-1}}{\pi\Gamma(b)\Gamma(c-b)}\int\limits_{0}^{1}
\frac{t^{b-1}(1-t)^{c-b-1}\sin\left\{a\arctan\left(\frac{ty\sin(\pi/a)}{1+ty\cos(\pi/a)}\right)\right\}}
{\left(1+2ty\cos(\pi/a)+t^2y^2\right)^{a/2}}dt.
\]
In particular, for $a=2$ we obtain:
\[
{_{2}F_1}(2,b;c;-z)=\frac{4b}{\pi{c}}\int\limits_{0}^{\infty}\frac{y^{2}{_{3}F_2}(2,(b+1)/2,(b+2)/2;(c+1)/2,(c+2)/2;-y^2)dy}{y^2+z^2},
\]
where we have used
\[
\int\limits_{0}^{1}\frac{t^{b}(1-t)^{c-b-1}}
{\left(1+t^2y^2\right)^{2}}dt=
\frac{\Gamma(b+1)\Gamma(c-b)}{\Gamma(c+1)}
{_{3}F_2}(2,(b+1)/2,(b+2)/2;(c+1)/2,(c+2)/2;-y^2).
\]

Using some known results and techniques representation
(\ref{eq:Frepr}) together with Lemma~\ref{lm:G-positive} and
Theorem~\ref{th:q1Fq-order} leads to a number of implications for
generalized hypergeometric function which we present in the
subsequent sections. All statements presented below  are believed
to be new.

\paragraph {3. Inequalities for ${_{q+1}F_q}$.} Many results of
\cite{KarpSitnik} are based on representation (\ref{eq:Frepr})
with non-negative $\rho$.  However, the inequality $\rho\geq{0}$
has only been proved in this reference for $b_k>a_k>0$,
$k=1,2,\ldots,q$. Theorem~\ref{th:q1Fq-order} combined with some
results of \cite{KarpPril} allow us to extend the results of
\cite{KarpSitnik}  to all values of $a_k$, $b_k$ satisfying
(\ref{eq:amajorb}). In particular, we get the following
statements.
\begin{theorem}\label{th:monoton}
Suppose $B\prec^W{A}$ and $\delta>0$. Then the function
\begin{equation}\label{eq:f-defined}
x\to \frac{{_{q+1}F_q}\left(\sigma,A+\delta;
B+\delta;-x\right)}{{_{q+1}F_q}\left(\sigma,A;B;-x\right)}
\end{equation}
is monotone decreasing on $(-1,\infty)$ if $\sigma>0$ and monotone
increasing if $\sigma<0$.
\end{theorem}
The proof of this result in \cite[Theorem~1]{KarpSitnik} is based
on representation (\ref{eq:Frepr}) with non-negative $\rho$ and so
it applies to our situation here if $B\prec^W{A}$ and $\psi>0$.
The claim is then extended by continuity to $\psi=0$.

Next, we obtain a lower bound.
\begin{theorem}\label{th:LowerGeneral}
Suppose $B\prec^W{A}$ and $\sigma>0$.  Then  for all $x>-1$ the
inequality
\begin{equation}\label{eq:LowerGeneral}
\frac{1}{\left(1+x\prod_{i=1}^{q}(a_i/b_i)\right)^{\sigma}}\leq{_{q+1}F_q}(\sigma,A;B;-x)
\end{equation}
holds true with equality only for $x=0$.
\end{theorem}
\textbf{Proof.} Consider the case $0<\sigma\leq{1}$ first. Then
according to Theorem~\ref{th:q1Fq-order} and
\cite[Theorem~12]{KarpPril} the condition $B\prec^W{A}$ implies
that the function $[{_{q+1}F_q}(\sigma,A;B;-x)]^{1/\sigma}$
belongs to $S_1$. Note that the condition
$\sigma\leq\min(a_1,\ldots,a_q)$ from Theorem~\ref{th:q1Fq-order}
is not required to make this conclusion.  It is immediate to check
that
\[
\frac{1}{1+x\prod_{i=1}^{q}(a_i/b_i)}
\]
is the Pad\'{e} approximation to
$[{_{q+1}F_q}(\sigma,(a_q);(b_q);-x)]^{1/\sigma}$ at $x=0$ of
order $[0/1]$.  This implies (\ref{eq:LowerGeneral}) for all
$x>-1$ by Stieltjes inequalities \cite[formulas (3),
(4)]{Gilewicz}.

Next, suppose that $\sigma>1$.  Then (\ref{eq:LowerGeneral}) can
be derived from Theorem~\ref{th:monoton} by repeating the proof of
\cite[Theorem~3]{KarpSitnik} word for word.~~$\square$

Inequality (\ref{eq:LowerGeneral}) was probably first obtained by
Luke in \cite{Luke} for $x>0$ and $b_k\geq{a_k}>0$.
Theorem~\ref{th:LowerGeneral}  extends his result to all $x>-1$
and parameters satisfying much weaker restrictions $B\prec^W{A}$.
An  extension of \cite[Theorem~4]{KarpSitnik} reads:
\begin{theorem}\label{th:UpperGeneral}
Suppose $B\prec^W{A}$ and $a_1,b_1>1$. Then for $x>0$ and
$0<\sigma\leq{1}$ the inequality
\begin{equation}\label{eq:UpperGeneral}
{_{q+1}F_q}(\sigma,A;B;-x)<\frac{1}{\left(1+x\prod_{i=1}^{q}[(a_i-1)/(b_i-1)]\right)^\sigma}
\end{equation}
holds.
\end{theorem}

In \cite{KarpSitnik1} Karp and Sitnik gave sufficient conditions
for absolute monotonicity of certain product differences of the
functions $_{q+1}F_q$.  This type of  absolute monotonicity
immediately implies log-convexity or log-concavity of
$\sigma\to_{q+1}F_q(\sigma,A;B;x)$ for $0<x<1$. Representation
(\ref{eq:Frepr}) allows for extension of log-convexity to $x<0$
under the restriction $B\prec^W{A}$.

\begin{theorem}\label{th:log-conv}
Suppose $B\prec^W{A}$.  Then the function
\[
\sigma\to{_{q+1}F_q}(\sigma,A;B;x)=:f(\sigma)
\]
is log-convex on $[0,\infty)$ for each $x<1$.
\end{theorem}
\textbf{Proof.}  Take $\sigma_2>\sigma_1\geq{0}$ and arbitrary
$\delta>0$.  The inequality
\[
f(\sigma_1+\delta)f(\sigma_2)\leq f(\sigma_1)f(\sigma_2+\delta)
\]
is equivalent to log-convexity for continuous functions (and is
stronger in general, see \cite[Chapter I.4]{MPF}), so it suffices
to prove this inequality.  Substituting (\ref{eq:Frepr}) for
$f(\sigma)$ we see that the above inequality is an instance of the
Chebyshev inequality \cite[Chapter IX, formula (1.1)]{MPF} if we
choose
\[
p(s)=\frac{\rho(s)}{(1-sx)^{\sigma_1}},~~~f(s)=\frac{1}{(1-sx)^{\sigma_2-\sigma_1}},~~~g(s)=\frac{1}{(1-sx)^{\delta}}.
\]
Indeed, $p(s)\geq{0}$ and both $f(s)$ and $g(s)$ are decreasing on
$(0,1)$ if $x<0$ and increasing if $0<x<1$. ~~$\square$

Some comments are in order here.  Using a completely different
approach Karp and Sitnik proved Theorem~\ref{th:log-conv} in
\cite{KarpSitnik1} for $0<x<1$ under the following conditions on
parameters:
\begin{equation}\label{eq:symmetric-chain}
\frac{e_q(b_1,\ldots,b_q)}{e_q(a_1,\ldots,a_q)}\geq
\frac{e_{q-1}(b_1,\ldots,b_q)}{e_{q-1}(a_1,\ldots,a_q)}\geq\cdots\geq
\frac{e_1(b_1,\ldots,b_q)}{e_1(a_1,\ldots,a_q)}\geq{1}
\end{equation}
where
\[
e_k(x_1,\ldots,x_q)=\sum\limits_{1\leq{j_1}<{j_2}\cdots<{j_k}\leq{q}}x_{j_1}x_{j_2}\cdots{x_{j_k}}
\]
is $k$-th elementary symmetric polynomial.  It is curious to
compare the conditions (\ref{eq:amajorb}) and
(\ref{eq:symmetric-chain}).  The essential part of this comparison
was done by Issai Schur in 1923.  More precisely, we have
\begin{lemma}\label{lm:schur}
Suppose $B\prec^W{A}$.  Then \emph{(\ref{eq:symmetric-chain})}
holds.
\end{lemma}
\textbf{Proof.}   According to \cite[3.A.8]{MOA} $B\prec^W{A}$
implies that $\phi(A)\leq\phi(B)$ if and only if $\phi(x)$ is
Schur-concave and increasing in each variable. Inequalities
(\ref{eq:symmetric-chain}) can alternatively be written as
\[
\frac{e_{k}(a_1,\ldots,a_q)}{e_{k-1}(a_1,\ldots,a_q)}\leq\frac{e_{k}(b_1,\ldots,b_q)}{e_{k-1}(b_1,\ldots,b_q)},~~k=1,2,\ldots,q.
\]
So we should choose
\[
\phi_k(x_1,\ldots,x_q)=\frac{e_{k}(x_1,\ldots,x_q)}{e_{k-1}(x_1,\ldots,x_q)},~~k=1,2,\ldots,q.
\]
Schur-concavity of these functions has been proved by Schur (1923)
- see \cite[3.F.3]{MOA}.  It is left to show that $\phi_k$ is
increasing in each variable.  Due to symmetry we can take $x_1$ to
be variable thinking of $x_2,\ldots,x_q$ as being fixed.  Using
the definition of elementary symmetric polynomials we see that for
$k\geq{2}$
\[
\phi_k(x_1,\ldots,x_q)=\frac{x_1e_{k-1}(x_2,\ldots,x_q)+e_{k}(x_2,\ldots,x_q)}{x_1e_{k-2}(x_2,\ldots,x_q)+e_{k-1}(x_2,\ldots,x_q)}.
\]
So taking derivative with respect to  $x_1$ we obtain
($e_m=e_m(x_2,\ldots,x_q)$ for brevity):
\[
\frac{\partial\phi_k(x_1,\ldots,x_q)}{\partial{x_1}}=\frac{e_{k-1}(x_1e_{k-2}+e_{k-1})-e_{k-2}(x_1e_{k-1}+e_{k})}
{[x_1e_{k-2}+e_{k-1}]^2}=\frac{e_{k-1}^2-e_{k}e_{k-2}}{[x_1e_{k-2}+e_{k-1}]^2}\geq{0}.
\]
Non-negativity holds by Newton's inequalities.~~$\square$

\textbf{Remark.} Since the reverse implication in
Lemma~\ref{lm:schur} is clearly not true, we see that the
log-convexity of ${_{q+1}F_q}(\sigma,A;B;x)$ in $\sigma$ holds for
$x<0$ under the conditions $B\prec^W{A}$ and for $0\leq{x}<1$
under weaker conditions (\ref{eq:symmetric-chain}).  Numerical
experiments show that the log-convexity indeed  does not hold for
$x<0$ under conditions (\ref{eq:symmetric-chain}) if we violate
$B\prec^W{A}$.

\paragraph{4. Pad\'{e} approximation to ${_{q+1}F_q}$.}
Theorem~\ref{th:q1Fq-order} together with
\cite[Theorem~3]{KarpPril} imply that for $B\prec^W{A}$ and
$0<\sigma\leq{1}$ the function
\[
z\to{_{q+1}F_q}(\sigma,A;B;-z):={_{q+1}F_q}(-z)
\]
belongs to the Stieltjes cone $S_1$ with the representing measure
$\rho$ supported on $[0,1]$.   This fact has a number of
consequences for the Pad\'{e} table of ${_{q+1}F_q}$.  Before
stating them we give an explicit expression for the density which
follows directly from Theorem~\ref{th:FGrep}.

\begin{theorem}
Suppose $\Re\left(\sum_{i=1}^{q}(b_i-a_i)\right)+1>\Re(\sigma)$
and $\Re(a_i)>0$, $i=1,\ldots,q$.  Then
\begin{equation}\label{eq:q1F1Stieltjes}
{_{q+1}F_q}\left(\left.\!\!\begin{array}{c}\sigma,A\\
B\end{array}\right|-z\!\right)
=\!\!\int\limits_{0}^{1}\frac{\rho_1(s)ds}{1+sz},
\end{equation}
with
\begin{equation*}
\rho_1(s)=\frac{1}{\Gamma(\sigma)}\left(\prod\limits_{i=1}^{q}\frac{\Gamma(b_i)}{\Gamma(a_i)}\right)
\frac{1}{s}G^{q+1,0}_{q+1,q+1}\left(s\,\,\vline\begin{array}{c}1,B\\\sigma,A\end{array}\!\!\right).
\end{equation*}
\end{theorem}
\textbf{Proof.}  Write
\[
{_{q+1}F_q}\left(\left.\!\!\begin{array}{c}\sigma,A\\
B\end{array}\right|-z\!\right)
={_{q+2}F_{q+1}}\left(\left.\!\!\begin{array}{c}1,\sigma,A\\
1,B\end{array}\right|-z\!\right)
\]
and apply Theorem~\ref{th:FGrep}.~~$\square$

Representation (\ref{eq:q1F1Stieltjes}) leads to:
\begin{theorem}\label{th:Pade-normal}
Suppose  $B\prec^W{A}$ and $0<\sigma\leq{1}$. Then for all integer
$m,n\geq{0}$ the Pad\'{e} approximant $[m/n]$ to ${_{q+1}F_q}(-z)$
at $z=0$ is normal.
\end{theorem}
\textbf{Proof.} Follows from representation
(\ref{eq:q1F1Stieltjes}) by
\cite[Theorem~4.2.3]{CPVWJ}.~~$\square$

\textbf{Remark.}  Let us remind the reader that a Pad\'{e}
approximant is called normal if it occupies precisely one entry in
the Pad\'{e} table.

\begin{theorem}\label{th:Pade-convergence}
Suppose  $B\prec^W{A}$ and $0<\sigma\leq{1}$.  Then the Pad\'{e}
approximants $[m+j/m]$, $j\geq{-1}$, converge to ${_{q+1}F_q}(-z)$
uniformly on every compact subset of $\C\setminus(-\infty,-1]$ as
$m\to\infty$.
\end{theorem}
\textbf{Proof.} Follows from representation
(\ref{eq:q1F1Stieltjes}) by \cite[Theorem~5.4.2]{BG}.~~$\square$

\begin{theorem}\label{th:Pade-denominators}
Suppose $B\prec^W{A}$, $\psi>0$ and $0<\sigma\leq{1}$. Then the
Pad\'{e} approximants $[m+j/m]$, $j\geq{-1}$, to ${_{q+1}F_q}(-z)$
have the form
\[
\frac{P^{[m+j/m]}(z)}{Q^{[m+j/m]}(z)}=\frac{P^{[m+j/m]}(z)}{(-z)^m\pi^j_m(-1/z)},
\]
where $\pi_m^j(s)$  are polynomials orthogonal  with respect to
the following inner product:
\[
\int\limits_{0}^{1}\pi_m^j(s)\pi_n^j(s)s^jG^{q+1,0}_{q+1,q+1}\left(s\,\,\vline\begin{array}{c}1,B\\\sigma,A\end{array}\!\!\right)ds=
\mathrm{const}\times\delta_{mn}.
\]
The numerator polynomials $P^{[m+j/m]}(z)$ are found from
\[
{_{q+1}F_q}(\sigma,A;B;-z)Q^{[m+j/m]}(z)-P^{[m+j/m]}(z)=O(z^{2m+j+1}),~~z\to{0}.
\]
\end{theorem}
\textbf{Proof.} Follows from representation
(\ref{eq:q1F1Stieltjes}), Lemma~\ref{lm:G-positive} and
\cite[Chapter 5, formula (3.21)]{BG}.~~$\square$

\paragraph{5. Mapping properties of ${_{q+1}F_q}$.} There is a vast
literature dedicated to the mapping properties of the Gauss
hypergeometric function ${_2F_1}$.  However, the mapping
properties of the functions ${_{q+1}F_q}(z)$ and $z{_{q+1}F_q}(z)$
for $q\geq{2}$ have been only considered by a few authors
\cite{OS,Ponnusamy}.  A combination of \cite[Theorem~13,
Remark~7]{KarpPril} with Theorem~\ref{th:q1Fq-order} immediately
yields

\begin{theorem}\label{th:pFq-sectormapping}
Suppose $B\prec^W{A}$ and $\sigma\geq{1}$. Then the function
${_{q+1}F_q}(\sigma,A;B;-z)$ maps the sector
$0<\arg(z)<\pi/\sigma$ into the lower half-plane $\Im(z)<0$.
\end{theorem}

Here we only demonstrate the direct consequences of
Theorem~\ref{th:q1Fq-order} when it is combined with the results
of Thale \cite{Thale} and Wirths \cite{Wirths}.
\begin{theorem}\label{th:Univ-starlike}
Suppose $B\prec^W{A}$ and $0<\sigma\leq{1}$.  Then the functions
\[
z\to{_{q+1}F_q}(\sigma,A;B;z)~~\text{and}~~z\to
z{_{q+1}F_q}(\sigma,A;B;z)
\]
are univalent in the half-plane $\Re(z)<1$. The second function is
also starlike in the disk $|z|<r^*$, where
\[
r^*=\sqrt{13\sqrt{13}-46}\approx{0,934}.
\]
\end{theorem}
The proof of the first claim follows from representation
(\ref{eq:Frepr}) combined with \cite[Theorems~2.1, 2.2]{Thale} or
\cite[Satz~2.2]{Wirths}. The second claim follows from
\cite[Satz~2.4]{Wirths}.~~$\square$

\textbf{Remark.} The constant $r^{*}$ above looks different from
the (much more cumbersome) constant given in \cite{Wirths} but a
simple calculation shows that they are equal.

\begin{theorem}\label{th:univ2}
Suppose $B\prec^W{A}$ and $0<\sigma\leq{2}$.  Then the function
\[
z\to z{_{q+1}F_q}(\sigma,A;B;z)
\]
is univalent in the disk $|z|<r_s:=\sqrt{\sqrt{32}-5}\approx
0.81$.
\end{theorem}
The claim follows from representation (\ref{eq:Frepr}) combined
with \cite[Satz~3.2]{Wirths}.

\bigskip

\paragraph{6. Acknowledgements.} We thank  Sergei
Sitnik (Voronezh Institute of the Ministry of Internal Affairs of
the Russian Federation) for useful discussions. We acknowledge the
financial support of the Russian Basic Research Fund (grant
11-01-00038-a) and the Far Eastern Branch of the Russian Academy
of Sciences.

\end{document}